\documentclass{elsart}
\usepackage{amsmath, amssymb, latexsym}

\newtheorem{theorem}{Theorem}

\newcommand{\N}{\ensuremath{\mathbb{N}}}
\newcommand{\R}{\ensuremath{\mathbb{R}}}

\def \a {\alpha}

\def \e {\varepsilon}
\def \d {\delta}
\def \f {\phi}
\def \l {\lambda}
\def \m {\mu}

\def \< {\langle}
\def \> {\rangle}
\def \ra {\rightarrow}

\def \span {{\rm span}}

\newtheorem{definition}{Definition}

\newcommand{\C}{\mathbb{C}}

\newcommand{\adj}{^\ast}

\newcommand{\inv}{^{-1}}

\newcommand{\LL}{\mathcal{L}}

\newenvironment{proof}[1][Proof]{\setlength{\parskip}{0pt plus 0.4pt}
\normalfont\upshape
\trivlist
\item[\hskip\labelsep{\itshape#1.}\ignorespaces]}{\qed\endtrivlist\setlength{\parskip}{0pt plus 1pt}}

\begin{document}
\begin{frontmatter}

\title{On the relation of closed forms and Trotter's product formula}
\author{M\'at\'e Matolcsi}
\address{ Alfr\'ed R\'enyi Institute of Mathematics,
                                     Hungarian Academy of Sciences
                                     POB 127
                                     H-1364 Budapest, Hungary
                                     Tel: (+361) 483-8302
                                     Fax: (+361) 483-8333, 
e-mail: matomate@renyi.hu   }

\begin{abstract}
The aim of this paper is to give a characterization in Hilbert spaces of the generators of 
$C_0$-semigroups associated with closed, sectorial forms in terms of the convergence of a 
generalized Trotter's product formula. In the course of the proof of the main result we also
 present a similarity result which can be of independent interest: for any unbounded generator
 $A$ of a $C_0$-semigroup $e^{tA}$ it is possible to introduce an equivalent scalar product 
on the space, such that $e^{tA}$ becomes non-quasi-contractive with respect to the new scalar 
product.
\end{abstract}

\end{frontmatter} 

\section{Introduction}\label{sec1}
Let $H$ denote a complex Hilbert space, $A$ the generator of a $C_0$-semigroup $e^{tA}$ on $H$, 
and $P$ a bounded projection. The convergence of the generalized Trotter's product formula
\begin{eqnarray}\label{pro}  
\lim_{n \rightarrow \infty}(e^{\frac{t}{n}A}P)^n
\end{eqnarray}
was first studied (as a corollary of the main Theorem) in \cite{Ka2}. It was shown that \eqref{pro} converges strongly for all $t>0$ 
whenever $-A$ is a non-negative self-adjoint operator, and $P$ is an orthogonal projection. More 
generally, the Addendum in \cite{Ka2} implies that the same result is true whenever $-A$ is 
associated with a closed sectorial form (where the vertex of the sector is allowed to be any 
real number $\omega$, i.e. 
$S_{\f ,\omega}:=\{z\in \C : -\f <{\mathrm{arg}} \ (z-\omega)< 
\f\}$, $\omega\in\R$, $\f\in (0,\frac{\pi}{2})$).

The convergence of \eqref{pro} was then studied in more general settings in \cite{AB} and 
\cite{MS} where further convergence results, motivating examples, and some counterexamples were 
given. The main result of this paper is to prove the converse of Kato's result, i.e. that the 
strong convergence of \eqref{pro} for all orthogonal projections $P$, in fact, characterizes 
generators $A$ such that $-A$ is associated with a closed sectorial form. To be more precise 
we recall the following result (see \cite{Ka2} Addendum, and \cite{MS} Theorem 1 and 4.):
\begin{theorem}\label{thm1}
Let $A$ be the generator of a $C_0$-semigroup $e^{tA}$ on a Hilbert 
space $H$. Consider the following statements:

(i) $A$ is bounded.

(ii) $-A$ is associated with a densely-defined, closed, sectorial form $a$ on $H$. 

(iii) The formula $(e^{\frac{t}{n}A}P)^nx$ converges for all projections $P\in \LL (H)$, 
and all $x\in H$ and $t>0$.

(iv) The formula $(e^{\frac{t}{n}A}P)^nx$ converges for all {\it{orthogonal}} projections
 $P\in \LL (H)$, and all $x\in H$ and $t>0$.

The following implications hold: (i) $\Rightarrow$ (iii)  and (ii) $\Rightarrow$ (iv).
\end{theorem}

We will show in Section \ref{sec2} that the converse implications also hold. In the course of 
the proof we will need an auxillary result, given in Theorem \ref{thm2}, which can be regarded 
as a complement of \cite{Ch}, and is of independent interest. Namely, we show that whenever the generator 
$A$ of the semigroup $e^{tA}$ is unbounded, it is possible to introduce an equivalent scalar product 
$( \ , \ )_0$ on $H$ such that $e^{tA}$ is non-quasi-contractive with respect to $( \ , \ )_0$.

\section{Main Result}\label{sec2}

In order to prove our main result [Theorem \ref{thm3}], first we need to characterize the class 
of generators $A$ on $H$, such that the $C_0$-semigroup $e^{tA}$ is quasi-contractive for 
{\it{every}} equivalent scalar product $( \ , \ )_0$ on $H$. The characterization is provided by 
\begin{theorem}\label{thm2}
Let $A$ be the generator of a $C_0$-semigroup $e^{tA}$ on a Hilbert space $H$.
 The following are equivalent:

(i) A is bounded.

(ii) The semigorup $e^{tA}$  is quasi-contractive for every equivalent scalar product 
$( \ , \ )_0$ on $H$.

(iii)  For every equivalent scalar product $( \ , \ )_0$ on $H$ there exists 
$K_0\in \R$ such that for every vector $x\in D(A)$, $(x,x)_0=1$ 
implies ${\mathrm{Re}} \ (Ax,x)_0\le K_0$.

\begin{proof}
The implications (ii) $\Leftrightarrow$ (iii) are consequences of the Lumer-Phillips theorem 
(see e.g. \cite{RN}, Proposition 3.23.). The implications (i) $\Rightarrow$ (ii) and (i) 
$\Rightarrow$ (iii) are trivial. It remains to prove (iii) $\Rightarrow$ (i). We will need the 
following 
\begin{definition}\label{def1}
Let $T\in \LL (H)$ be an injective operator, and $x\in H$, $\|x\|=1$, and $0<\d \le 1$. 
We say that $x$ is a $\d$-quasi-eigenvector of $T$ if   
\begin{eqnarray}
\d \le \frac{|(x,Tx)|}{\|Tx\|}\le 1
\end{eqnarray}
\end{definition}
Note, that a 1-quasi-eigenvector is, in fact, an eigenvector of $T$. 

Now, let $0<\d <1$ be fixed. We prove the implication (iii) $\Rightarrow$ (i) by contradiction. 
Assume, therefore, that $A\notin \LL (H)$, and also, by rescaling, that $A^{-1}=:T\in\LL (H)$. 
Assume, furthermore, that a sequence $(h_n)\subset H$ is given with the following properties:

(a) $\|h_n\|=1$ for all $n\ge 1$.

(b) $\{h_k, \  Th_k\} \perp  \{h_j, \ Th_j\}$ for all $k \not= j$.

(c) $\lim_{n\to\infty}\|Th_n\|=0$   

(d) For every $n\ge 1$ the vector $h_n$ is {\it{not}} a $\d$-quasi-eginvector of $T$.

We construct an equivalent scalar product $( \ , \ )_0$ on $H$ with the help of the sequence 
$h_n$.

Let $H_n=\span\{h_n, \ Th_n\}$. Note, that $H_n$ is 2-dimensional because $h_n$ is not an 
eigenvector of $T$.

Let $Th_n=c_{1,n}h_n+c_{2,n}h_n^{\perp}$, where $\|h_n^{\perp}\|=1$.  
Note that 
$$\frac{|c_{1,n}|^2}{|c_{1,n}|^2+|c_{2,n}|^2}<\d^2  \ \  {\mathrm{and}} \  \  
\frac{|c_{2,n}|^2}{|c_{1,n}|^2+|c_{2,n}|^2}>1-\d^2$$ 
Hence, 
$$\frac{|c_{1,n}|}{|c_{2,n}|}<\frac{\d}{\sqrt{1-\d^2}} \ \  {\mathrm{and}} \ \  
\frac{|c_{2,n}|}{\|Th_n\|}>\sqrt{1-\d^2}$$ 

Define $Q_n\in\LL (H_n)$ by 
\begin{eqnarray*} 
Q_nh_n:=h_n+\overline{L}_nh_n^{\perp}\\
Q_nh_n^{\perp}:=L_nh_n+(|L_n|^2+1)h_n^{\perp}
\end{eqnarray*}
where $|L_n|=2\frac{\d}{\sqrt{1-\d^2}}$ and $\overline{{L}_n{c}_{2,n}}>0$ for 
all $n\ge 1$. It is clear that $Q_n=Q_n\adj \ge 0$, $Q_n\inv\in \LL(H_n)$, 
and $\|Q_n\|_{H_n}\le K$,  $\|Q_n\inv\|_{H_n}\le K$ for some universal constant $K$ (not 
depending on $n$). Define $Q\in \LL (H)$ by 
\begin{eqnarray*}
Q:=Q_1\oplus Q_2\oplus\dots\bigoplus I_{(H_1\oplus H_2\oplus\dots)^{\perp}}
\end{eqnarray*}
It is easy to see that $Q$ is well-defined, $Q\in\LL (H)$, $Q=Q\adj\ge 0$, and 
$Q\inv\in\LL (H)$. This means that $Q$ defines an equivalent scalar product on $H$ by 
$(x,y)_0:=(x,Qy)$. 

Now, let $x_n:=\frac{Th_n}{\|Th_n\|}$. Then 
\begin{eqnarray*}
{\mathrm{Re}} \ (Ax_n,x_n)_0=\frac{1}{\|Th_n\|^2}{\mathrm{Re}} \ (h_n,QTh_n)=
\frac{1}{\|Th_n\|^2}{\mathrm{Re}} \ (h_n,c_{1,n}h_n+c_{2,n}L_nh_n)=\\
=\frac{1}{\|Th_n\|^2}({\mathrm{Re}} \ c_{1,n}+\overline{c_{2,n}L_n})\ge
\frac{1}{\|Th_n\|^2}\frac{\d}{\sqrt{1-\d^2}}|c_{2,n}|\ge \frac{1}{\|Th_n\|}\d\ra +\infty
\end{eqnarray*}
Let $y_n:=\frac{x_n}{\|x_n\|}_0$. Then $ {\mathrm{Re}} \ (Ay_n,y_n)_0\ra +\infty$ still holds
 due to the equivalence of the scalar products $( \ , \ )$ and $( \ , \ )_0$. 

In order to complete the proof of the theorem it remains to construct the sequence $h_n$ with 
the required properties. The construction is carried out in several steps. 

{\it{Step 1.}} We construct an orthonormal sequence $(e_n)\subset H$, such that\\
$\lim_{n\to\infty}\|Te_n\|=0$. 

Take the polar decomposition $T=UT_1$ of $T$, where $U$ is unitary and $T_1=T_1\adj\ge 0$. It
 is clear from the spectral theorem that there exists an orthonormal sequence $(e_n)\subset H$ 
such that 
$\lim_{n\to\infty}\|T_1e_n\|=0$ (otherwise $T_1$ and $T$ would be invertible, contrary
 to our assumption). Note, also, that $\|T_1e_n\|= \|Te_n\|$ for all $n\in\N$, therefore 
$\lim_{n\to\infty}\|Te_n\|=0$ as required.

{\it{Step 2.}} We construct an orthonormal sequence $(f_n)\subset H$ such that \\ 
$\lim_{n\to\infty}\|Tf_n\|=0$ and $f_{n+1}\perp \{f_1,Tf_1,...f_n,Tf_n\}$. 

We obtain the sequence $(f_n)$ by induction, with the help of the sequence $(e_n)$. Take an 
index $i_1$ such that  $\|Te_{i_1}\|\le 1$, and let $f_1:=e_{i_1}$. Assume now that 
$f_1, f_2,\dots ,f_n$ are already given such that 
$$\|f_j\|=1, \ \ f_j\perp \{f_k,Tf_k\}, \  \ \|Tf_j\|\le\frac{1}{\sqrt{j}}$$ 
and $f_j\in\span \{e_1,e_2,\dots ,e_{l_n}\}$, for all $1\le j,k\le n$, $k<j$, and $l_n$ is an 
index depending on $n$ only.

Let $H_n:=\span\{Tf_1,Tf_2,\dots ,Tf_n\}$. Take indices $j_1,\dots j_{n+1}$ such that $j_k>l_n$ 
and  $\|Te_{j_k}\|\le\frac{1}{n+1}$ for all $1\le k\le n+1$. The subspace $H_n$ is at most 
$n$-dimensional, therefore there exists a non-trivial linear combination 
$$f_{n+1}:=\sum_{k=1}^{n+1}\l_ke_{j_k}$$ 
such that $\|f_{n+1}\|=1$ and $f_{n+1}\perp H_n$. 

It is clear, by construction, that $f_{n+1}\perp \{f_1,Tf_1,...f_n,Tf_n\}$. Furthermore, 
$$\|Tf_{n+1}\|\le\frac{1}{n+1}\sum_{k=1}^{n+1}|\l_k|\le
\sqrt{\frac{\sum_{k=1}^{n+1}|\l_k|^2}{n+1}}=\frac{1}{\sqrt{n+1}}$$ 

{\it{Step 3.}}  We construct an orthonormal sequence $(g_n)\subset H$ such that\\ 
$\lim_{n\to\infty}\|Tg_n\|=0$ and $\{g_j, Tg_j\}\perp \{g_k,Tg_k\}$ for all $j\not= k$.

We obtain the sequence $(g_n)$ by induction, with the help of the sequence $(f_n)$.

Let $g_1=f_1$. 
Assume now that $g_1, g_2,\dots ,g_n$ are already given such that 
$$\|g_j\|=1, \ \  \{g_j,Tg_j\}\perp \{g_k,Tg_k\}, \ \ \|Tg_j\|\le\frac{1}{\sqrt{2j-1}}$$ 
and $g_j\in\span \{f_1,f_2,\dots ,f_{b_n}\}$, for all $1\le j\not= k\le n$, and $b_n$ is an 
index depending on $n$ only.

Let $G_n:=\span\{g_1, Tg_1, g_2, Tg_2,\dots ,g_n, Tg_n\}$. Take indices $m_1,\dots m_{2n+1}$ 
such that $m_k>b_n$ and  $\|Tf_{m_k}\|\le\frac{1}{2n+1}$ for all $1\le k\le 2n+1$. The subspace 
$G_n$ is at most $2n$-dimensional, therefore there exists a non-trivial linear combination 
$$g_{n+1}:=\sum_{k=1}^{2n+1}\m_kf_{m_k}$$ 
such that $\|g_{n+1}\|=1$ and $Tg_{n+1}\perp G_n$.    

It is clear, by construction, that $\{g_{n+1},Tg_{n+1}\}\perp \{g_1,Tg_1,...g_n,Tg_n\}$. 
Furthermore, 
$$\|Tg_{n+1}\|\le\frac{1}{2n+1}\sum_{k=1}^{2n+1}|\m_k|\le
\sqrt{\frac{\sum_{k=1}^{2n+1}|\m_k|^2}{2n+1}}=\frac{1}{\sqrt{2(n+1)-1}}$$

{\it{Step 4.}} We construct the orthonormal sequence $(h_n)$ with the properties stated at the 
beginning of the proof. 

We obtain the sequence $(h_n)$ by induction, with the help of the sequence $(g_n)$.

Take an index $r_1$ such that $\|Tg_{r_1}\|\le \frac{\d^2}{10}\|Tg_1\|$. Let 
$$h_1:=\frac{\d}{2}g_1+\sqrt{1-\frac{\d^2}{4}}g_{r_1}$$ 
We need to prove that $h_1$ is not a $\d$-quasi-eigenvector of $T$. It is clear that 
$$1\ge \|Th_1\|\ge \left( \frac{\d}{2}-\frac{\d^2}{10}\sqrt{1-\frac{\d^2}{4}}\right)\|Tg_1\|$$ 
Also, 
$$|(h_1,Th_1)|= |(g_1,Tg_1)+(g_{r_1},Tg_{r_1})|\le \left(\frac{\d^2}{4}+(1-\frac{\d^2}{4})
\frac{\d^2}{10}\right)\|Tg_1\|$$
Combining these two ineqalities a simple calculation shows that 
$\frac{|(h_1,Th_1)|}{\|Th_1\|}< \d$, as required.

Assume now that vectors $h_1, \dots ,h_n$ are already given, such that $h_j$ is not a 
$\d$-quasi-eigenvector of $T$, 
$$\|h_j\|=1, \ \   \{h_j,Th_j\}\perp \{h_k,Th_k\}, \ \  \|Th_j\|\le\frac{1}{\sqrt{j}}$$ 
and $h_j\in\span \{g_1,g_2,\dots ,g_{a_n}\}$, for all $1\le j\not= k\le n$, and $a_n$ is an 
index depending on $n$ only. Take indices $p_1$, $p_2$, such that $p_1,p_2>a_n$ 
and $\|Tg_{p_1}\|\le\frac{1}{\sqrt{n+1}}$, and  $\|Tg_{p_2}\|\le \frac{\d^2}{10}\|Tg_{p_1}\|$. 
Let 
$$h_{n+1}:=\frac{\d}{2}g_{p_1}+\sqrt{1-\frac{\d^2}{4}}g_{p_2}$$ 
It is clear that $\|Th_{n+1}\|\le \frac{1}{\sqrt{n+1}}$, and it can be shown as above 
that $h_{n+1}$ is not a $\d$-quasi-eigenvector of $T$. Hence, the sequence $(h_n)$ satisfies 
all requirements, and the proof is complete.
\end{proof}
\end{theorem}

The author conjectures that a result corresponding to Theorem \ref{thm2} holds
also in Banach spaces. Namely, whenever $A$ is not bounded it should be
possible to introduce an equivalent norm on the space such that $e^{tA}$ is
not quasi-contractive with respect to the new norm. This problem, however,
remains open.   

Now we present the main result of the paper.

\begin{theorem}\label{thm3}
Let $A$ be the generator of a $C_0$-semigroup $e^{tA}$ on a Hilbert space $H$. 
Consider the following statements. 

(i) A is bounded.

(ii) $-A$ is associated with a densely-defined, closed, sectorial form $a$ on $H$. 

(iii) The formula $(e^{\frac{t}{n}A}P)^nx$ converges for all projections $P\in \LL (H)$, and 
all $x\in H$ and $t>0$.

(iv) The formula $(e^{\frac{t}{n}A}P)^nx$ converges for all {\it{orthogonal}} projections 
$P\in \LL (H)$, and all $x\in H$ and $t>0$.

The following implications hold: (i) $\Leftrightarrow$ (iii), (ii) $\Leftrightarrow$ (iv).

\begin{proof}
The implication  (i) $\Rightarrow$ (iii) was proved in \cite{MS}, while the implication\\ 
(ii) $\Rightarrow$ (iv) is a consequence of \cite{Ka2}, Addendum (see also \cite{MS}, Theorem 4).

We prove the implication (iii) $\Rightarrow$ (i) by contradiction. 

Assume first that the semigroup $e^{tA}$ is not quasi-contractive. By the Lumer-Phillips 
theorem this is equivalent to the fact that the numerical range of $A$ is not contained in any 
left half-plane. 

We construct an element $g\in H$ such that $\|g\|=1$, and 
$$\lim_{n\ra\infty}(e^{\frac{1}{n}A}P_g)^ng$$ 
does not exist, where $P_g$ denotes the one-dimensional projection onto the subspace 
spanned by $g$. The vector $g$ will be given as 
$$g:=\frac{\lim_{k\ra\infty}g_k}{\|\lim_{k\ra\infty}g_k\|}$$ 
where $(g_k)$ denotes a convergent sequence in $H$ to be constructed in the sequel.

Let $g_1\in D(A)$, such that $\|g_1\|=1$. First, we show that    
$$\lim_{n\ra\infty}(e^{\frac{1}{n}A}P_{g_1})^ng_1=e^{(Ag_1,g_1)}g_1$$ 
Indeed, 
$$(e^{\frac{1}{n}A}P_{g_1})^ng_1=e^{\frac{1}{n}A}(P_{g_1}e^{\frac{1}{n}A}P_{g_1})^{n-1}g_1=
e^{\frac{1}{n}A}(P_{g_1}e^{\frac{1}{n}A}P_{g_1}g_1, \ g_1)^{n-1}g_1$$
and 
$$\lim_{n\ra\infty}(P_{g_1}e^{\frac{1}{n}A}P_{g_1}g_1, \ g_1)^{n-1}=e^{(Ag_1, \ g_1)}$$
because 
$$\lim_{n\ra\infty}\frac{(P_{g_1}e^{\frac{1}{n}A}P_{g_1}g_1, \ g_1)-1}{1/n}=\lim_{n\ra\infty}
\left( \frac{(e^{\frac{1}{n}A}-I)g_1}{1/n}, \ g_1\right)=(Ag_1, \ g_1)$$

Now, choose $g_1$ such that ${\mathrm{Re}} \ (Ag_1,g_1)\ge 1$ holds also.

Let $\e >0$ be fixed. Take an index $n_1$ so large that 
$$\| \left(e^{\frac{1}{n_1}A}P_{g_1}\right)^{n_1}g_1-e^{(Ag_1,g_1)}g_1\|<\e$$ 
It is clear from standard continuity arguments that there exists a $\d_1 >0$, such that for 
all $h\in B(g_1,\d_1)$ we have 
$$\| \left(e^{\frac{1}{n_1}A}P_{\frac{h}{\|h\|}}\right)^{n_1}\frac{h}{\|h\|}-
e^{(Ag_1,g_1)}g_1\|<2\e$$ 
Without loss of generality we can assume that $\d_1<\frac{1}{2}$. 

Now assume, that vectors $g_1, g_2, \dots , g_k$, and positive numbers 
$\d_1,\d_2,\dots , \d_k$, and indices $n_1,n_2,\dots ,n_k$  are already given with  
the properties that: 
$$g_j\in D(A), \ \  {\mathrm{Re}} \ (Ag_j,g_j)\ge j$$ 
and  
$$\| \left(e^{\frac{1}{n_j}A}P_{\frac{h}{\|h\|}}\right)^{n_j}\frac{h}{\|h\|}-
e^{(A\frac{g_j}{\|g_j\|},\frac{g_j}{\|g_j\|})}\frac{g_j}{\|g_j\|}\|<2\e$$ 
for all $1\le j\le k$ and all $h\in B(g_j,\d_j)$. Assume, furthermore, that 
$$\|g_{j+1}-g_j\|< {\mathrm{min}} \ \left\{\frac{\d_1}{2^j}, \frac{\d_2}{2^{j-1}}, \dots  
\frac{\d_j}{2}\right\}$$ 
for all $1\le j\le k-1$. 

The numerical range of $A$ is not bounded from the right, hence there exists a vector 
$f\in D(A)$ such that 
$$\|f\|<  {\mathrm{min}} \ \left\{\frac{1}{\|Ag_k\|}, \frac{\d_1}{2^k}, 
\frac{\d_2}{2^{k-1}}, \dots  \frac{\d_k}{2}\right\}$$ 
and 
${\mathrm{Re}}\ \ (Af,f)\ge 2$. 
Let $f_k:=e^{i\a}f$ with suitable $\a$ such that  ${\mathrm{Re}}\ \ (Af_k,g_k)\ge 0$. Let 
$$g_{k+1}:=g_k+f_k$$ 
Then  
\begin{eqnarray*}
{\mathrm{Re}}\ \ (Ag_{k+1},g_{k+1})=
{\mathrm{Re}}\ \ (Ag_{k},g_{k}) + {\mathrm{Re}}\ \ (Ag_{k},f_{k}) +\\ 
+{\mathrm{Re}}\ \ (Af_{k},g_{k}) + {\mathrm{Re}}\ \ (Af_{k},f_{k})
\ge k+(-1)+0+2=k+1 
\end{eqnarray*}
Furthermore, we have 
$$\lim_{n\to\infty} (e^{\frac{1}{n}A}P_{ \frac{g_{k+1}}{\|g_{k+1}\|}  })^{n} 
\frac{g_{k+1}}{\|g_{k+1}\|}=e^{(A\frac{g_{k+1}}{\|g_{k+1}\|},\frac{g_{k+1}}{\|g_{k+1}\|})}
\frac{g_{k+1}}{\|g_{k+1}\|}$$ 
Take an index $n_{k+1}$ so large that $n_{k+1}>n_k$ and 
$$\| \left(e^{\frac{1}{n_{k+1}}A}P_{ \frac{g_{k+1}}{\|g_{k+1}\|}   }\right)^{n_{k+1}} 
\frac{g_{k+1}}{\|g_{k+1}\|}  -e^{(A \frac{g_{k+1}}{\|g_{k+1}\|}  , 
\frac{g_{k+1}}{\|g_{k+1}\|}  )} \frac{g_{k+1}}{\|g_{k+1}\|}  \|<\e$$ 
It is clear from standard continuity arguments that there exists a $\d_{k+1} >0$, such that for 
all $h\in B(g_{k+1},\d_{k+1})$ we have 
$$\|\left(e^{\frac{1}{n_{k+1}}A}P_{\frac{h}{\|h\|}}\right)^{n_{k+1}}\frac{h}{\|h\|}-
e^{(A \frac{g_{k+1}}{\|g_{k+1}\|}  , \frac{g_{k+1}}{\|g_{k+1}\|}  )} 
\frac{g_{k+1}}{\|g_{k+1}\|}  \|<2\e$$

It is clear, by construction, that the sequence $g_k$ converges in $H$. Let 
$$h:=\lim_{k\to\infty}g_k \ \ {\mathrm{and}} \ \  g:= \frac{h}{\|h\|}$$ 
Recall, that $\|g_1\|=1$ and $\d_1<\frac{1}{2}$, therefore $\frac{1}{2}<\|g_k\|<\frac{3}{2}$ 
for all $k\ge 1$. It is also clear, by construction, that  $h\in B(g_k,\d_k)$ for all $k\ge 1$. 
Hence, for all $k\ge 1$ we have  
$$\| \left(e^{\frac{1}{n_k}A}P_{g}\right)^{n_k}g-e^{(A\frac{g_k}{\|g_k\|},
\frac{g_k}{\|g_k\|})}\frac{g_k}{\|g_k\|}\|<2\e$$ 
Notice, that 
$$\|e^{(A\frac{g_k}{\|g_k\|},\frac{g_k}{\|g_k\|})}\frac{g_k}{\|g_k\|}\|=
e^{\frac{1}{\|g_k\|^2}{\mathrm{Re}} \ (Ag_k,g_k)}> e^{\frac{1}{4}k}$$ 
This means that (the norm of) the sequence 
$ (e^{\frac{1}{n}A}P_{g})^{n}g$ does not converge.

Now, assume only that $A\notin \LL (H)$. Introduce, by Theorem \ref{thm2}, an equivalent scalar 
product $(x,y)_0:=(x,Qy)$ on $H$, such that the semigroup $e^{tA}$ is not quasi-contractive 
with respect to  $( \ ,  \ )_0$. Take an orthogonal projection $P_g$ (with respect to the 
scalar product $( \ ,  \ )_0 \ $),  such that  $ (e^{\frac{1}{n}A}P_{g})^{n}g$ does not 
converge. Then, $P_g$ is a bounded (possibly non-orthogonal) projection with respect to the 
original scalar product $(  \ ,  \ )$, such that 
$ (e^{\frac{1}{n}A}P_{g})^{n}g$ does not converge. This proves the implication 
(iii) $\Rightarrow$ (i). 

The implication (iv)  $\Rightarrow$ (ii) is also proved by contradiction.

Asume, that the numerical range of $A$ is not contained in any sector 
$$\Sigma_{\f ,\omega}:=\{z\in \C : \frac{\pi}{2}+\f <{\mathrm{arg}} \ (z-\omega)< 
\frac{3}{2}\pi-\f\}$$ 
with $\omega\in\R$, $\f\in (0,\frac{\pi}{2})$. There are two cases to consider. 

If the semigroup $e^{tA}$ is not quasi-contractive , then, by the arguments above, 
there exists a vector $g\in H$, such that $\|g\|=1$ and  $ (e^{\frac{1}{n}A}P_{g})^{n}g$ 
does not converge.

If the semigroup $e^{tA}$ is quasi-contractive then, by rescaling, we can assume that 
${\mathrm{Re}} \ (Ax,x)\le -1$ for all $x\in D(A)$, $\|x\|=1$.  

We construct an element $g\in H$ such that $\|g\|=1$, and 
$\lim_{n\ra\infty}(e^{\frac{1}{n}A}P_g)^ng$ does not exist, where $P_g$ denotes the 
one-dimensional projection onto the subspace spanned by $g$. The vector $g$ will be given as 
$$g:=\frac{\lim_{k\ra\infty}g_k}{\|\lim_{k\ra\infty}g_k\|}$$ 
where $(g_k)$ denotes a convergent sequence in $H$ to be constructed in the sequel.

Take an arbitrary vector $g_1\in D(A)$, $\|g_1\|=1$. 
Let $(Ag_1,g_1)=:a_1+b_1i$. We know that 
$$\lim_{n\ra\infty}\left(e^{\frac{1}{n}A}P_{g_1}\right)^ng_1=e^{(Ag_1,g_1)}g_1$$ 
Let $\e>0$, and $\rho >0$ be fixed. Take an index $n_1$ so large that 
$$\| \left(e^{\frac{1}{n_1}A}P_{g_1}\right)^{n_1}g_1-e^{(Ag_1,g_1)}g_1\|<\e$$ 
It is clear from standard continuity arguments that there exists a $\d_1 >0$, such that for 
all $h\in B(g_1,\d_1)$ we have 
$$\| \left(e^{\frac{1}{n_1}A}P_{\frac{h}{\|h\|}}\right)^{n_1}\frac{h}{\|h\|}-
e^{(Ag_1,g_1)}g_1\|<2\e$$ 
Without loss of generality we can assume that $\d_1<\frac{1}{2}$. 

Now assume, that vectors $g_1, g_2, \dots , g_k$, and positive numbers 
$\d_1,\d_2,\dots , \d_k$,\\ 
real numbers $\e_1 , \e_2 , \dots ,\e_k$, and indices $n_1,n_2,\dots ,n_k$  are already given 
with the following properties: for all $1\le j\le k$ we have  $|\e_j|<\rho$,  
$$g_j\in D(A), \ \ (A\frac{g_j}{\|g_j\|}, \frac{g_j}{\|g_j\|}  )= a_j+(\e_j+b_1+(j-1)\pi )i$$ 
(note that $\e_1=0$),
where $a_1-1<a_j\le -1$, and  
$$\| \left(e^{\frac{1}{n_j}A}P_{\frac{h}{\|h\|}}\right)^{n_j}\frac{h}{\|h\|}-
e^{(A\frac{g_j}{\|g_j\|},\frac{g_j}{\|g_j\|})}\frac{g_j}{\|g_j\|}\|<2\e$$ 
for all $h\in B(g_j,\d_j)$. Assume, furthermore, that 
$$\|g_{j+1}-g_j\|< {\mathrm{min}} \ \left\{\frac{\d_1}{2^j}, \frac{\d_2}{2^{j-1}}, \dots  
\frac{\d_j}{2}\right\}$$ for all $1\le j\le k-1$. 

Now, we construct the vector $g_{k+1}$. The numerical range of $A$ is not contained in any 
sector, therefore there exists a sequence $(x_j)\subset D(A)$ such that, 
$\lim_{j\to\infty}\|x_j\|=0$ and 
$${\mathrm{Im}} \ \frac{(Ax_j,x_j)}{ \|g_k\|^2   }=
\pi \ \ {\mathrm{and}} \ \  \frac{{\mathrm{Re}} \ (Ax_j,x_j)}{\|g_k\|^2}<\frac{a_k-(a_1-1)}{2}$$ 

Take $y_j:=e^{i\a_j}x_j$ with suitable $\a_j$ such that $(Ay_j,g_k)\ge 0$ real. Then 
\begin{eqnarray*}
\frac{(A(g_k+y_j), g_k+y_j)}{\|g_k\|^2}=\frac{(Ag_k,g_k)}{\|g_k\|^2}+
\frac{(Ag_k,y_j)}{\|g_k\|^2}+\\
+\frac{(Ay_j,g_k)}{\|g_k\|^2}+\frac{(Ay_j,y_j)}{\|g_k\|^2}=:c_j+d_ji 
\end{eqnarray*}
The real part $c_j$ of this expression satisfies 
$$c_j>(a_1-1)+(\frac{a_k-(a_1-1)}{2})-\frac{|(Ag_k,y_j)|}{\|g_k\|^2}$$ 
for all $j\ge 1$. For the imaginary part $d_j$, we have 
$$\lim_{j\to\infty}d_j=\e_k+b_1+k\pi$$ 
This means that for large $j$ we have $\|y_j\|< {\mathrm{min}} \ \{\frac{\d_1}{2^k}, 
\frac{\d_2}{2^{k-1}}, \dots  \frac{\d_k}{2}\}$, and 
$$\frac{{\mathrm{Re}}(A(g_k+y_j),g_k+y_j)}{\|g_k+y_j\|^2}>a_1-1$$ 
and  
$$\frac{{\mathrm{Im}}(A(g_k+y_j),g_k+y_j)}{\|g_k+y_j\|^2}=\e_{k+1}+b_1+k\pi$$ 
where $|\e_{k+1}|<\rho$. Take such an index $j$, and define 
$$g_{k+1}:=g_k+y_j$$ 
Again, standard continuity arguments show that there exist a positive number $\d_{k+1}$ and an 
index $n_{k+1}$ such that 
$$\| \left(e^{\frac{1}{n_{k+1}}A}P_{\frac{h}{\|h\|}}\right)^{n_{k+1}}\frac{h}{\|h\|}-
e^{(A\frac{g_{k+1}}{\|g_{k+1}\|},\frac{g_{k+1}}{\|g_{k+1}\|})}\frac{g_{k+1}}{\|g_{k+1}\|}\|<2\e$$ 
for all $h\in B(g_{k+1},\d_{k+1})$. 

It is clear, by construction, that the sequence $g_k$ converges.  Let 
$$h:=\lim_{k\to\infty}g_k \ \ {\mathrm{ and}} \ \ g:= \frac{h}{\|h\|}$$ 
Recall, that $\|g_1\|=1$ and $\d_1<\frac{1}{2}$, therefore $\frac{1}{2}<\|g_k\|<\frac{3}{2}$ 
for all $k\ge 1$. It is also clear, by construction, that  $h\in B(g_k,\d_k)$ for all $k\ge 1$. 
Hence, for all $k\ge 1$ we have  
$$\| \left(e^{\frac{1}{n_k}A}P_{g}\right)^{n_k}g-e^{(A\frac{g_k}{\|g_k\|},
\frac{g_k}{\|g_k\|})}\frac{g_k}{\|g_k\|}\|<2\e$$

Notice, furthermore that 
\begin{eqnarray*}
\|e^{(A\frac{g_{2k+1}}{\|g_{2k+1}\|},\frac{g_{2k+1}}{\|g_{2k+1}\|})}
\frac{g_{2k+1}}{\|g_{2k+1}\|}-e^{(A\frac{g_{2k}}{\|g_{2k}\|},\frac{g_{2k}}{\|g_{2k}\|})}
\frac{g_{2k}}{\|g_{2k}\|}\|=\\ 
\|e^{a_{2k+1}}e^{(\e_{2k+1}+b_1+2k\pi)i}\frac{g_{2k+1}}{\|g_{2k+1}\|}-
e^{a_{2k}}e^{(\e_{2k}+b_1+(2k-1)\pi)i}\frac{g_{2k}}{\|g_{2k}\|}\|\ge\\ 
\|e^{a_{2k+1}+b_1i}g_1-e^{a_{2k}+(b_1-\pi)i}g_1\|-
\|e^{a_{2k+1}+b_1i}(e^{\e_{2k+1}}\frac{g_{2k+1}}{\|g_{2k+1}\|}-g_1)\|-\\
-\|e^{a_{2k}+(b_1-\pi)i}(e^{\e_{2k}}\frac{g_{2k}}{\|g_{2k}\|}-g_1)\|\ge\\
 2e^{a_1-1}-\|e^{a_{2k+1}+b_1i}(e^{\e_{2k+1}}\frac{g_{2k+1}}{\|g_{2k+1}\|}-g_1)\|-
\|e^{a_{2k}+(b_1-\pi )i}(e^{\e_{2k}}\frac{g_{2k}}{\|g_{2k}\|}-g_1)\| 
\end{eqnarray*}
We can now choose the values of $\e ,\d_1, \rho$ so small that 
\begin{eqnarray*}
\|e^{a_{2k+1}+b_1i}(e^{\e_{2k+1}}\frac{g_{2k+1}}{\|g_{2k+1}\|}-g_1)\|+
\|e^{a_{2k}+(b_1-\pi)i}(e^{\e_{2k}}\frac{g_{2k}}{\|g_{2k}\|}-g_1)\|\le e^{a_1-1} 
\end{eqnarray*}
and $5\e \le e^{a_1-1}$ 

Then we have 
$$\| (e^{\frac{1}{n_{2k+1}}A}P_{g})^{n_{2k+1}}g- (e^{\frac{1}{n_{2k}}A}P_{g})^{n_{2k}}g\|\ge \e$$ 
Therefore the sequence $(e^{\frac{1}{n}A}P_{g})^{n}g$ does not converge, and the proof is 
complete. 

\end{proof}
\end{theorem}

\end{document}